\documentclass[11pt]{article}

\title{The exponentially weighted average forecaster in geodesic spaces of non-positive curvature}
\author{Quentin Paris\footnote{HSE University, Faculty of Computer Science, Moscow, Russia. This work has been funded by the  Russian Academic Excellence Project '5-100'. Email:\url{qparis@hse.ru}}}
\date{}

\usepackage{amsthm,amssymb,amsmath,eucal}
\usepackage{color}
\usepackage{enumitem}
\usepackage{url}
\usepackage{bbold}
\usepackage{cases}
\usepackage{float}

\usepackage[top=2.5cm,bottom=2.5cm,left=2.5cm,right=2.5cm]{geometry}

\usepackage{natbib}
\numberwithin{equation}{section}
\usepackage{hyperref}
\hypersetup{
  colorlinks = true,
  urlcolor = blue, 
  linkcolor = blue,
  citecolor = blue}

\newtheorem{thm}{Theorem}[section]

\newtheorem{lem}[thm]{Lemma}
\newtheorem{cor}[thm]{Corollary}
\newtheorem{defi}[thm]{Definition}
\newtheorem{rem}[thm]{Remark}
\newtheorem{exm}[thm]{Example}

\newcommand{\prob}{\mathbb{P}}
\newcommand{\esp}{\mathbb{E}}

\newcommand{\R}{\mathbb R}

\newcommand{\curv}{\mathrm{curv}}

\begin{document}

\maketitle
\begin{abstract}
 This paper addresses the problem of prediction with expert advice for outcomes in a geodesic space with non-positive curvature in the sense of Alexandrov. Via geometric considerations, and in particular the notion of barycenters, we extend to this setting the definition and analysis of the classical exponentially weighted average forecaster. We also adapt the principle of online to batch conversion to this setting. We shortly discuss the application of these results in the context of aggregation and for the problem of barycenter estimation.\vspace{0.1cm}
\end{abstract}

\section{Introduction}
The problem of prediction with expert advice~\citep{CesLug06} is a by now standard model of online learning. Traditionally studied for outcomes taking values in a vector space, less seems to be known when the outcome space is a more general metric space. This paper partly addresses the problem by focusing on the case of NPC spaces, i.e., geodesic metric spaces with non-positive curvature in the sense of Alexandrov. 

The class of NPC spaces includes many metric spaces of particular interest in the data sciences. Apart from Hilbert spaces, interesting examples are hyperbolic spaces~\citep{NiKi17}, the space of real symmetric positive-definite matrices with Log-Euclidean~\citep{AFPA2007} or Log-Cholesky~\citep{Lin2019} Riemannian metrics and more generally all complete and simply connected Riemannian manifolds with non-positive sectional curvature. Another example is the BHV space of phylogenetic trees~\citep{bhv}. Finally, this class includes many other non-smooth objects such as metric trees or euclidean buildings, as well as images, products and gluings of spaces of non-positive curvature. For more details, we refer the reader to~\citet{BriHaf99,BurBurIva01,Stu03,AleKapPet19_invitation,AleKapPet19} and the references therein. 

Via geometric considerations, and in particular the notion of barycenters, we extend to this setting the definition and analysis of the classical exponentially weighted average (EWA) forecaster. In particular, our results rely on the notions of geodesic convexity, barycenters and a generalized version of Jensen's inequality instead of usual convexity, classical averages and the standard version of Jensen's inequality respectively. 

Another focus of the paper is the adaptation of the online-to-batch conversion principle, relying once again on the notion of barycenters. In the context of supervised learning with response variable taking values in an NPC space, we show that, combined with the generalized EWA forecaster, these results can be used to derive bounds for the aggregation of predictors. As a second application, we derive new bounds for the estimation of barycenters in NPC spaces.

\subsection{Prediction with expert advice}
\label{subsec:pea}
We first recall the protocol of prediction with expert advice. We refer the reader to~\citet{CesLug06} for more details. 

The values of some unknown sequence $z_1,z_2,\dots$ are revealed one by one in some outcome set $\mathcal Z$. At each round $t\ge 1$, a forecaster is to provide a prediction $\hat m_t\in \mathcal M$ for $z_t$, before it is revealed. To construct his prediction, the forecaster is given (in addition to past outcomes $(z_s)_{1\le s\le t-1}$) access to expert predictions $(m_{\theta,t})\in \mathcal M$ indexed by a set of experts $\Theta$. Hence, at each time step $t\ge 1$, the prediction $\hat m_t$ of the forecaster is formally constructed as 
\begin{equation}
\label{eq:predstrat}
\hat m_t=m_t((m_{\theta,t})_{\theta\in\Theta},(((m_{\theta,s})_{\theta\in\Theta},z_s))_{1\le s\le t-1}),
\end{equation}
for some prediction strategy $m_t:\mathcal M^{\Theta}\times (\mathcal M^{\Theta}\times \mathcal Z)^{t-1}\to \mathcal M$. Depending on the problem at hand, the set $\mathcal M$, referred to as the decision set, can be different from the outcome set $\mathcal Z$. Given a loss function 
\[\ell:\mathcal M\times \mathcal Z\to \R,\] 
fixed in advance, the cumulative performances of the forecaster and expert $\theta\in \Theta$ over a period of time $T\ge 1$ are respectively measured by
\[\hat L_T:=\sum_{t=1}^T\ell(\hat m_t,z_t)\quad\mbox{and}\quad L_{\theta,T}:=\sum_{t=1}^T\ell(m_{\theta,t},z_t).\]
The goal of the forecaster is to construct a prediction strategy minimizing the regret
\[R_T:=\hat L_T-\min_{\theta\in\Theta}L_{\theta,T},\]
uniformly over all outcome sequences and all expert advice. 

\subsection{The exponentially weighted average forecaster}
\label{subsec:ewaclassic}
In this learning scenario, it is classical to consider the decision set $\mathcal M$ to be a convex subset of $\R^p$, for some $p\ge 1$. In this case, a popular prediction strategy is the EWA forecaster defined as follows.

First, let $\pi$ be a prior distribution over the expert set $\Theta$ and $(\beta_t)_{t\ge 1}$ be a sequence of positive tuning parameters. Then, for $t\ge 1$, the EWA forecaster is  
\begin{equation}
\label{eq:ewaclassic}
\hat m_t=\int_{\Theta}m_{\theta,t}\,\mathrm{d}\pi_t(\theta),
\end{equation}
where 
\begin{equation}
\label{eq:ewaclassicw}
\mathrm{d}\pi_t(\theta)=\frac{e^{-\beta_t L_{\theta,t-1}}\,\mathrm{d}\pi(\theta)}{\int_{\Theta}e^{-\beta_t L_{\vartheta,t-1}}\,\mathrm{d}\pi(\vartheta)},
\end{equation}
and where $L_{\theta,0}=0$ by convention. This popular forecaster  naturally emphasizes the role of experts that perform well over time and its analysis is simplified by the convenient properties of the exponential function. 

Since its introduction by~\citet{Vov90} and~\citet{LitWar94}, the EWA forecaster has been analyzed from many perspectives~\citep{CesFreHausHelmSchWar97,CesLug99,Ces99}. The use of exponential weights has also found many successful statistical applications, in particular in the context of aggregation~\citep{Yan04,LeuBar06,Cat07,DalTsy07,DalTsy08,DalTsy09,JudRigTsy08,Alq08,Aud09,DalTsy12a,DalTsy12b}.   

Among the many results available in the literature, the following PAC-bayesian result is now standard and has the advantage of allowing for an arbitrary expert set $\Theta$ (we recommend~\citealt[Chapter 2]{Ger11}, for a comprehensive review of similar results). Given a probability measure $q$ on $\Theta$, we denote $D(q\|\pi)$ the Kulback-Leibler divergence between $q$ and $\pi$, i.e., 
\[D(q\|\pi):=\int_{\Theta}\ln\left(\frac{{\mathrm d}q}{\mathrm{d}\pi}\right){\mathrm d}q,\]
if $q\ll\pi$ and $D(q\|\pi):=+\infty$ otherwise.

\begin{thm}
\label{thm:pac}
Suppose there exists $\beta>0$ such that, for all $z\in \mathcal Z$, the function $e^{-\beta\ell(.,z)}:\mathcal M\to \R$ is concave. Suppose $\hat m_t$ is constructed as in \eqref{eq:ewaclassic} with $\beta_t=\beta$ at each time step. Then, uniformly over the outcome sequence and the expert advice,
\[\hat L_T\le \inf_q\left\{\int_{\Theta}L_{\theta,T}\,{\mathrm d} q(\theta)+\frac{D(q\|\pi)}{\beta}\right\},\]
where the inf runs over all probability measures $q$ on $\Theta$. In particular if the expert set $\Theta$ is finite, with $K$ elements, and $\pi$ is the uniform prior over $\Theta$, we get the regret bound
\[R_T\le \frac{\ln K}{\beta}.\]
\end{thm}

Even thought classical, the proof is shortly outlined in Appendix \ref{sec:proofs} for its informative content. In the classical setting considered here, the linear structure of the decision space $\mathcal M$ is important for several reasons. First, it allows to make sense of the integral defining $\hat m_t$ in \eqref{eq:ewaclassic}. Second, it is essential to the usual notion of convexity invoked in the statement of Theorem \ref{thm:pac} and in Jensen's inequality used in the proof. As it turns out, a simple adaptation of the construction of $\hat m$ given in \eqref{eq:ewaclassic}, using the notion of barycenters, makes sense in an abstract metric space and reduces to the familiar EWA forecaster in the euclidean setting. Our goal in this paper will be to show in particular that the analysis of this generalized EWA forecaster can be carried out in a similar way in NPC spaces. 

\subsection{Organisation of the paper}
Section \ref{sec:geom} presents some simple tools from metric geometry. Subsection \ref{subsec:ewa} presents and analyses the performance of a generalized EWA forecaster. Paragraph \ref{subsec:otob} proves online-to-batch conversion results. We shortly discuss the application of these results in the context of aggregation, in subsection \ref{subsec:learn}, and in the context of barycenter estimation in subsection \ref{subsec:estimbary}. Proofs are reported to Appendix \ref{sec:proofs}.

\section{Geometric preliminaries}
\label{sec:geom}
 In this section, we mention some standard facts from metric geometry for reference and prove some useful results. For more details on metric geometry, we refer the reader to~\citet{BriHaf99,BurBurIva01,AleKapPet19_invitation,AleKapPet19} and the references therein.
 
\subsection{Geodesic spaces}
Let $(M,d)$ be a metric space. We call path in $M$ a continuous map $\gamma:I\to M$ defined on an interval $I\subset \R$. A path $\gamma:[0,1]\to M$ will be called a geodesic if, for all $0\le s\le t\le 1$,
\begin{equation}
\label{eq:scalegeod}
d(\gamma(s),\gamma(t))=(t-s)d(\gamma(0),\gamma(1)),
\end{equation}

\begin{defi}[Geodesic space]
A metric space $(M,d)$ is called geodesic if, for every $x,y\in M$, there exists a geodesic $\gamma:[0,1]\to M$ connecting $x$ to $y$, i.e., such that $\gamma(0)=x$ and $\gamma(1)=y$. 
\end{defi}

Any (convex subset of a) normed vector space $(V,\|.\|)$ is a geodesic space. For instance, given any $x,y\in V$, the path $\gamma(t)=(1-t)x+ty$ defines a geodesic\footnote{Note however that it needs not be the unique geodesic. One checks that it is unique if the norm $\|.\|$ is strictly convex, e.g., if $V$ is a Hilbert space with inner product norm.} connecting $x$ to $y$ according to characterization \eqref{eq:scalegeod}. Other fundamental examples of geodesic spaces are provided by complete and connected Riemannian manifolds $M$, with Riemannian metric $g$, equipped with the Riemannian distance 
\[d(x,y)=\inf_{\gamma}\int_{0}^{1}\sqrt{g_{\gamma(t)}(\gamma'(t),\gamma'(t))}\,\mathrm{d}t,\]
where the infimum is taken over all paths $\gamma:[0,1]\to M$ connecting $x$ to $y$. Finally, the class of geodesic spaces admits many other examples of objects which cannot be described as Riemannian manifolds because of there singularities or infinite dimensional nature. Examples of such spaces are metric graphs or euclidean buildings~\citep[Chapter 3]{BurBurIva01} or the $2$-Wasserstein space over (for example) a euclidean space~\citep[Section 7.2]{AmbGigSav08}. 
 
\subsection{NPC spaces}

\begin{defi}[NPC spaces]
\label{def:npc}
Let $(M,d)$ be a geodesic space. We say that $\curv(M)\le 0$ if, for every $p\in M$, every $x,y\in M$ and every geodesic $\gamma:[0,1]\to M$ connecting $x$ to $y$, we have 
\begin{equation}
\label{eq:npc}    
d^2(p,\gamma(t))\le (1-t)d^2(p,x)+td^2(p,y)-t(1-t)d^2(x,y).
\end{equation}
\end{defi}
A geodesic spaces $(M,d)$ satisfying $\curv(M)\le 0$ is often termed an NPC\footnote{NPC stands for non-positive curvature. NPC spaces are also called Hadamard space or CAT(0) space.} space.  Next are some examples. Much more examples can be found in the literature mentioned in the beginning of this section.
\begin{exm}[Hilbert spaces]
    In the context where $(M,d)$ is a Hilbert space, equipped with its natural metric, inequality \eqref{eq:npc} is an identity. A Banach space is known to satisfies \eqref{eq:npc} iff it is a Hilbert space~\citep[Proposition 3.5]{Stu03}.
\end{exm}
\begin{exm}[Riemmannian manifolds]
    When $(M,d)$ is a complete and simply-connected Riemmannian manifold with non-positive sectional curvature, equipped with its Riemannian distance, then $\curv(M)\le 0$. Non-positive sectional curvature is known to imply the comparison inequality \eqref{eq:npc} locally (i.e., for $p,x,y$ close enough) and the simple-connectedness allows for this property to hold globally (i.e., for all $p,x,y$). The standard example of a Riemmanian manifold satisfying the conditions of Definition \ref{def:npc} is the hyperbolic plane.  One another example of interest in a vast number of applications is the space $M=S^{+}_d(\R)$ of real $d\times d$ symmetric positive definite matrices. Several works have studied the Riemannian structure of $S^{+}_d(\R)$. In this direction, \citet{AFPA2007} and \citet{Lin2019} introduce respectively the Log-Euclidean and Log-Cholesky Riemannian scalar products, both leading to a Riemannian distance on $S^{+}_d(\R)$ satisfying \eqref{eq:npc}.
    \end{exm}
    \begin{exm}[$L^2$ spaces]
    \label{exm:l2spaces}
    Consider two metric spaces $(\mathcal X,d_{\mathcal X})$ and $(\mathcal Y,d_{\mathcal Y})$ and a Borel probability measure $\mu$ on $\mathcal X$. Define $M$ as the set of all Borel measurable functions $f:\mathcal X\to \mathcal Y$ such that
    \[\int_{\mathcal X}d^2_{\mathcal Y}(f(x),y)\,\mathrm{d}\mu(x)<+\infty,\]
    for all $y\in \mathcal Y$. When $(\mathcal Y,d_{\mathcal Y})$ is a geodesic space satisfying $\curv(\mathcal Y)\le 0$, the set $(M,d)$, equipped with metric $d$ defined by
    \[d^2(f,g):=\int_{\mathcal X} d^2_{\mathcal Y}(f(x),g(x))\,\mathrm{d}\mu(x),\]
    is geodesic and satisfies $\curv(M)\le 0$~\citep[Proposition 3.10]{Stu03}. This example will be of interest for some applications considered later in the paper. 
    \end{exm}

 \subsection{Geodesic convexity}

Let $(M,d)$ be geodesic. A function $f:M\to \R$ is called geodesically convex if, for every $x,y\in M$ and every geodesic $\gamma:[0,1]\to M$ connecting $x$ to $y$, the function $f\circ\gamma:[0,1]\to \R$ is convex. We'll say that a function $f:M\to \R$ is geodesically concave if $-f$ is geodesically convex. Given $\alpha\in\R$, a function $f:M\to \R$ is called geodesically $\alpha$-convex if, for every $x,y\in M$ and every geodesic $\gamma:[0,1]\to M$ connecting $x$ to $y$, the function 
\[t\in[0,1]\mapsto f(\gamma(t))-\frac{\alpha}{2} d^2(\gamma(0),\gamma(1))t^2, \]
is convex. These definitions reduce to the usual notions of convexity in the context of euclidean spaces. These definitions also allow to provide an alternative characterization of NPC spaces. Indeed, it follows from the Definition \ref{def:npc} that ${\rm curv}(M)\le 0$ iff, for all $p\in M$, the function $d^2(.,p):M\to \R$ is geodesically $2$-convex. 

Next is a useful lemma proved in Appendix \ref{sec:proofs}. 

\begin{lem}
\label{lem:lcisec}
Let $(M,d)$ be a complete geodesic space and $f:M\to\R$ be a given function. 
\begin{enumerate}
\item Suppose that $e^{-\beta f}$ is geodesically concave for some $\beta>0$. Then, the function $f$ is geodesically convex.
\item Suppose $f$ is geodesically $\alpha$-convex and $L$-Lipchitz for some $\alpha,L>0$. Then, the function $e^{-\beta f}$ is geodesically concave for all 
$0\le \beta\le \frac{\alpha}{L^2}$.
\end{enumerate}
\end{lem}
To avoid additional technicality, all metric spaces considered in the paper will be considered complete. Using the triangular inequality, and the second statement above, we get the following straightforward corollary.

\begin{cor}
\label{cor:d2isbetaconcave}
Let $(M,d)$ be a bounded NPC space with diameter $\mathrm{diam}(M)>0$. Then, for all $p\in M$, the function $e^{-\beta d^2(.,p)}:M\to \R$ is geodesically concave for all \[0\le \beta\le \frac{1}{2\mathrm{diam}(M)^2}.\]
\end{cor}

\subsection{Barycenters and Jensen's inequality}
We now define barycenters which will play a central role in the next section.
\begin{defi}
Given a metric space $(M,d)$, define $\mathcal P_2(M)$ as the set of Borel probability measures $P$ on $M$ such that, for all $x\in M$,
\begin{equation}
\label{eq:vf}
\int_M d^2(x,y)\,\mathrm{d}P(y)<+\infty.
\end{equation}
A barycenter of $P\in \mathcal P_2(M)$ is any 
\begin{equation}
\label{xstarintrobary}
x^{*}\in\underset{x\in M}{\arg\min}\int_M d^2(x,y)\,\mathrm{d}P(y).
\end{equation}
\end{defi}
Barycenters\footnote{Barycenters are also called Fr\'echet means, intrinsic means, $2$-means or centers of mass.} provide a generalization\footnote{Strictly speaking, a barycenter generalizes the notion of mean value of a probability measure with a finite second moment. Indeed, if for instance $(M,d)=(\R^p,\|.-.\|_2)$ and provided $\int\|x\|^2_2\,\mathrm{d}P(x)<+\infty$, it is well known that 
$x^*=\int x\,\mathrm{d}P(x)$ is the unique minimizer of 
\[x\in M\mapsto\int \|x-y\|^2_2\,\mathrm{d}P(y).\]} of the notion of mean value when $M$ has no linear structure. While alternative notions of mean value in a metric space have been proposed, barycenters are usually favored for their simple interpretation and constructive definition as solution of an optimization problem. Solving numerically this minimization problem is, in itself, still an active field of investigation. However, the problem is rather well studied in the context of NPC spaces~\citep{Jos95,May98,AmbGigSav08,Bac14}.  The question of existence and uniqueness of barycenters is important and has been addressed in a number of settings. The following statement~\citep[Theorem 4.9]{Stu03} summarizes the results we'll need in the following. 

 \begin{thm}
 \label{thm:sturm}
 Let $(M,d)$ be a complete metric space. Then, the following statements are equivalent. 
 \begin{enumerate}
     \item $(M,d)$ is an NPC space.
     \item Any probability measure $P\in \mathcal P_2(M)$ has a unique barycenter $x^{*}$ and, for all $x\in M$, 
    \begin{equation}
        \label{eq:vistrurm}
    d^2(x,x^*)\le \int_M (d^2(x,y)-d^2(x^*,y))\,\mathrm{d}P(y).
    \end{equation}
 \end{enumerate} 
 \end{thm}
In the context of Euclidean spaces, note that inequality \eqref{eq:vistrurm} holds as an identity. Next is a generalized Jensen inequality \citep[Theorem 6.2]{Stu03}.

\begin{thm}
\label{thm:jensen}
Let $(M,d)$ be an NPC space. Let $f:M\to \R$ be (lower semi-continous and) geodesically convex and $P\in\mathcal P_2(M)$. Then (when the rhs is well defined),   
\[f(x^*)\le\int_M f(x)\,\mathrm{d}P(x),\]
where $x^*$ is the barycenter of $P$.
\end{thm}

\section{Results}
\label{sec:main}

\subsection{A generalized EWA forecaster in NPC spaces}
\label{subsec:ewa}
We are now in position to describe the generalized EWA forecaster. We adopt the same notation as in paragraphs \ref{subsec:pea} and \ref{subsec:ewaclassic}. From now on, we suppose that the decision set $\mathcal M$ is an NPC space with metric $d:\mathcal M\times\mathcal M\to \R_+$. 

Let $\pi$ be a prior distribution over the expert set $\Theta$ and $(\beta_t)_{t\ge 1}$ be a sequence of positive tuning parameters. Then, for $t\ge 1$, we set 
\begin{equation}
\label{eq:ewametric}
\hat m_t\in\underset{m\in\mathcal M}{\arg\min}\int_{\Theta}d^2(m,m_{\theta,t})\,\mathrm{d}\pi_t(\theta),
\end{equation}
where
\begin{equation}
\nonumber
\mathrm{d}\pi_t(\theta)=\frac{e^{-\beta_t L_{\theta,t-1}}\,\mathrm{d}\pi(\theta)}{\int_{\Theta}e^{-\beta_t L_{\vartheta,t-1}}\,\mathrm{d}\pi(\vartheta)},
\end{equation}
and where $L_{\theta,0}=0$ by convention. Note that while this definition adapts the original construction \eqref{eq:ewaclassic} to the metric setting, the weight update mechanism is kept as it is and still makes perfect sense. 

An important observation is that $\hat m_t$ can be alternatively written as
\begin{equation}
\label{eq:ewametric2}
\hat m_t\in\underset{m\in\mathcal M}{\arg\min}\int_{\mathcal M}d^2(m,x)\,\mathrm{d}P_t(x), 
\end{equation}
where $P_t$ is the pushforward of $\pi_t$ by the map $\theta\in\Theta\mapsto m_{\theta,t}\in\mathcal M$. In other words, $\hat m_t$ is the unique barycenter of probability distribution $P_t$. 
\begin{exm} While the above construction makes sense when the expert set $\Theta$ is an arbitrary measurable space, we translate it in the finite setting for illustration. Suppose the expert set is finite, say $\Theta=\{1,\dots,K\}$. Let $\pi=(\pi_1,\dots,\pi_K)$ be the prior distribution over $\Theta$ and denote, for all $1\le \theta\le K$,
\[\pi_{\theta,t}:=\frac{\pi_{\theta}e^{-\beta_t L_{\theta,t-1}}}{\sum_{\vartheta=1}^K\pi_{\vartheta}e^{-\beta_t L_{\vartheta,t-1}}}.\]
Then the generalized EWA forecaster in this setting is
\[\hat m_t\in\underset{m\in\mathcal M}{\arg\min}\sum_{\theta=1}^K\pi_{\theta,t}d^2(m,m_{\theta,t}).\]
In other words, $\hat m_t$ is the unique barycenter of probability measure $P_t$ on $\mathcal M$ defined by
\[P_t=\sum_{\theta=1}^K\pi_{\theta,t}\delta_{m_{\theta,t}},\]
where $\delta_m$ is the Dirac mass at $m$.
\end{exm}
We'll implicitly suppose, throughout the rest of the paper, that the map $\theta\in\Theta\mapsto m_{\theta,t}\in\mathcal M$ is measurable and that $P_t\in\mathcal P_2(M)$\footnote{In the euclidean setting, this amounts to assume that, \[\int_{\Theta}\|m_{\theta,t}\|^2_2\,\mathrm{d}\pi_t(\theta)<+\infty.\]} for all $t\ge 1$.  

The next result shows that the generalized EWA forecaster benefits from the exact same theoretical guarantees as in the euclidean setting displayed in Theorem \ref{thm:pac}.
\begin{thm}
\label{thm:ewametric1}
Suppose that $(\mathcal M,d)$ is an NPC space and that there exists $\beta>0$ such that, for all $z\in \mathcal Z$, the function $e^{-\beta\ell(.,z)}:\mathcal M\to \R$ is geodesically concave. Suppose $\hat m_t$ is constructed as in \eqref{eq:ewametric} with $\beta_t=\beta$ at each time step. Then, uniformly over the outcome sequence and the expert advice,
\[\hat L_T\le \inf_q\left\{\int_{\Theta}L_{\theta,T}\,{\mathrm d} q(\theta)+\frac{D(q\|\pi)}{\beta}\right\},\]
where the inf runs over all probability measures $q$ on $\Theta$. In particular, if the expert set $\Theta=\{1,\dots, K\}$ is finite and $\pi$ is the uniform prior over $\Theta$, we get
\[R_T\le \frac{\ln K}{\beta}.\]
\end{thm}

The next result is a direct consequence of Theorem \ref{thm:ewametric1} and Corollary \ref{cor:d2isbetaconcave}.

\begin{cor}
\label{cor:ewametric1}
Suppose that $(\mathcal M,d)$ is a bounded NPC space. Suppose that $\mathcal Z=\mathcal M$ and that $\ell(m,z)=d^2(m,z)$. Suppose finally that $\hat m_t$ is constructed as in \eqref{eq:ewametric} with $\beta_t=1/2\mathrm{diam}(\mathcal M)^2$ at each time step. Then, uniformly over the outcome sequence and the expert advice,
\[\hat L_T\le \inf_q\left\{\int_{\Theta}L_{\theta,T}\,{\mathrm d} q(\theta)+2\,\mathrm{diam}(\mathcal M)^2D(q\|\pi)\right\},\]
where the inf runs over all probability measures $q$ on $\Theta$. In particular if the expert set $\Theta$ is finite, with $K$ elements, and $\pi$ is the uniform prior over $\Theta$, we get
\[R_T\le 2\,\mathrm{diam}(\mathcal M)^2\ln K.\]
\end{cor}

To conclude, we mention that similar results, holding in situations where the loss $\ell:\mathcal M\times\mathcal Z\to\R$ is only geodesically convex in its first argument can be obtained along the exact same lines.

\subsection{Online-to-batch conversion}
\label{subsec:otob}

The principle of online-to-batch conversion is a classical and powerful way to exploit algorithms, developed for sequential prediction, in the context of statistical learning. This section mentions an adaptation of this well known procedure in the context of NPC spaces via the notion of barycenters and the use of Theorem \ref{thm:jensen}. 

Consider the following statistical learning problem. Suppose that $(\mathcal M,d)$ is an NPC space, that $\mathcal Z$ is a arbitrary measurable space and let $\ell:\mathcal M\times \mathcal Z\to \R$ be fixed. Suppose given a collection (or batch) $\{Z_{i}\}_{i=1}^n$ of independent and identically distributed $\mathcal Z$-valued random variables with same distribution as (and independent from) a generic random variable $Z$. Finally, given a subset $\Theta\subset \mathcal M$, consider the task of constructing $\theta_n\in \mathcal M$ (possibly not belonging to $\Theta$) based on $\{Z_{i}\}_{i=1}^n$ and such that
\[\esp[\ell(\theta_n,Z)]-\inf_{\theta\in\Theta}\esp[\ell(\theta,Z)],\]
is as small as possible.

To that aim, take the problem of sequential prediction with expert advice considered so far with constant expert advice, i.e., $m_{\theta,t}=\theta$ for all $\theta\in\Theta$ and all $t\ge 1$. In this simple setting, we can simplify the formal representation \eqref{eq:predstrat} of a prediction strategy $(\hat m_{t})_{t\ge 1}$ and consider that $\hat m_1=m_1$, for some constant $m_1\in \mathcal M$ independent of the outcome sequence, and that $\hat m_t=m_t(z_1,\dots,z_{t-1})$ for some function $m_t:\mathcal Z^{t-1}\to\mathcal M$ whenever $t\ge 2$.

Now, suppose given a prediction strategy $(\hat m_t)_{t\ge 1}$ such that, for all $T\ge 1$, there exists $B_T>0$ satisfying 
\begin{equation}
\label{eq:boundforotbc}
\sum_{t=1}^T\ell(\hat m_t,z_t)-\inf_{\theta\in\Theta}\sum_{t=1}^T\ell(\theta,z_t)\le B_T,
\end{equation}
uniformly over the outcome sequence $(z_1,\dots,z_T)\in\mathcal Z^T$. 

Then, coming back to the statistical learning problem, consider $\theta_n$ to be the unique barycenter of the (random) probability measure 
\[\frac{1}{n+1}\delta_{m_1}+\frac{1}{n+1}\sum_{i=2}^{n+1}\delta_{m_i(Z_1,\dots,Z_{i-1})},\]
on $\mathcal M$, i.e.,
\begin{equation}
    \label{eq:defthetan}
\theta_n\in\underset{m\in\mathcal M}{\arg\min}\left\{d^2(m,m_1)+\sum_{i=2}^{n+1} d^2(m,m_i(Z_1,\dots,Z_{i-1}))\right\}.
\end{equation}
Then we have the following result.

\begin{thm}
\label{thm:otbc}
Suppose that $(\mathcal M,d)$ is an NPC space. Suppose that \eqref{eq:boundforotbc} holds and that, for all $z\in\mathcal Z$, the function $\ell(.,z):\mathcal M\to \R$ is geodesically convex. 
Then, for all $n\ge 1$,
\[\esp[\ell(\theta_n,Z)]-\inf_{\theta\in\Theta}\esp[\ell(\theta,Z)]\le \frac{B_{n+1}}{n+1}.\]
\end{thm}
The next corollary is immediate by combining Theorem \ref{thm:otbc}, Theorem \ref{thm:ewametric1} and using statement $1$ in Lemma \ref{lem:lcisec}.
\begin{cor}
\label{cor:otbc}
Suppose that $(\mathcal M,d)$ is an NPC space. Suppose that $\Theta$ is a finite subset of $\mathcal M$. Suppose that there exists $\beta>0$ such that, for all $z\in\mathcal Z$, the function $e^{-\beta\ell(.,z)}:\mathcal M\to \R$ is geodesically concave. Let $\theta_n$ be as in \eqref{eq:defthetan} where $m_1$ is the barycenter of the uniform measure on $\Theta$ and, for $2\le i\le n+1$, $m_i(Z_1,\dots,Z_{i-1})$ is the barycenter of the (random) probability measure
\[\sum_{\theta\in\Theta}\pi_{\theta,i}\delta_{\theta},\quad\mbox{where}\quad\pi_{\theta,i}=\frac{e^{-\beta \sum_{j=1}^{i-1}\ell(\theta,Z_j)}}{\sum_{\vartheta\in\Theta}e^{-\beta \sum_{j=1}^{i-1}\ell(\vartheta,Z_j)}}.\]
Then, for all $n\ge 1$,
\[\esp[\ell(\theta_n,Z)]-\inf_{\theta\in\Theta}\esp[\ell(\theta,Z)]\le \frac{\ln |\Theta|}{\beta(n+1)}.\]
\end{cor}

The next results are in the same spirit as Theorem \ref{thm:otbc} and Corollary \ref{cor:otbc}. Instead of supposing that inequality \eqref{eq:boundforotbc} holds, assume that for every $T\ge 1$ there are positive numbers $B_T(q)>0$, indexed by probability distributions $q$ on $\Theta$, such that 
\begin{equation}
\label{eq:boundforotbc2}
\sum_{t=1}^T\ell(\hat m_t,z_t)\le\inf_q\left\{ \int_{\Theta}\sum_{t=1}^T\ell(\theta,z_t)\,\mathrm{d}q(\theta)+B_T(q)\right\},
\end{equation}
uniformly over the outcome sequence $(z_1,\dots,z_T)\in\mathcal Z^T$. Then we have the following result. Below we denote 
\[\theta^*\in\underset{\theta\in\Theta}{\arg\min}\,\esp[\ell(\theta,Z)].\]
\begin{thm}
\label{thm:otbc2}
Suppose that $(\mathcal M,d)$ is an NPC space. Suppose that \eqref{eq:boundforotbc2} holds and that, for all $z\in\mathcal Z$, the function $\ell(.,z):\mathcal M\to \R$ is geodesically convex. Then the following statements hold.
\begin{enumerate}
    \item For all $n\ge 1$,
\[\esp[\ell(\theta_n,Z)]\le\inf_q\left\{ \int_{\Theta}\esp[\ell(\theta,Z)]\,\mathrm{d}q(\theta)+\frac{B_{n+1}(q)}{n+1}\right\},\]
where the inf runs over all probability measures $q$ on $\Theta$.
\item If in addition there exists $\lambda>0$ such that, for all $z\in\mathcal Z$, $\ell(.,z):\mathcal M\to\R$ is $\lambda$-Lipschitz, then for all $n\ge 1$
\[\esp[\ell(\theta_n,Z)]-\inf_{\theta\in\Theta}\esp[\ell(\theta,Z)]\le\inf_q\left\{\lambda\int_{\Theta}d(\theta^*,\theta)\,\mathrm{d}q(\theta)+\frac{B_{n+1}(q)}{n+1}\right\},\]
where the inf runs over all probability measures $q$ on $\Theta$.
\end{enumerate}
\end{thm}

In particular, using the generalized EWA forecaster, we get the following.
\begin{cor}
\label{cor:otbc2}
Suppose that $(\mathcal M,d)$ is an NPC space. Let $\pi$ be a prior distribution over $\Theta$. Suppose that there exists $\beta>0$ such that, for all $z\in\mathcal Z$, the function $e^{-\beta\ell(.,z)}:\mathcal M\to \R$ is geodesically concave. Let $\theta_n$ be as in \eqref{eq:defthetan} where $m_1$ is the barycenter of $\pi$ and, for $2\le i\le n+1$, $m_i(Z_1,\dots,Z_{i-1})$ is the barycenter of the (random) probability measure $\pi_{i}$ defined by
\[\mathrm{d}\pi_{i}(\theta)=\frac{e^{-\beta \sum_{j=1}^{i-1}\ell(\theta,Z_j)}\mathrm{d}\pi(\theta)}{\int_{\Theta}e^{-\beta \sum_{j=1}^{i-1}\ell(\vartheta,Z_j)}\mathrm{d}\pi(\vartheta)}.\]
Then, the following statements hold. 
\begin{enumerate}
\item For all $n\ge 1$,
\[\esp[\ell(\theta_n,Z)]\le \inf_{q}\left\{\int_{\Theta}\esp[\ell(\theta,Z)]\,\mathrm{d}q(\theta)+\frac{D(q\|\pi)}{\beta(n+1)}\right\},\]
where the inf runs over all probability measures $q$ on $\Theta$.
\item If in addition there exists $\lambda>0$ such that, for all $z\in\mathcal Z$, $\ell(.,z):\mathcal M\to\R$ is $\lambda$-Lipschitz, then for all $n\ge 1$
\begin{equation}
\label{eq:otbc2}
    \esp[\ell(\theta_n,Z)]-\inf_{\theta\in\Theta}\esp[\ell(\theta,Z)]\le\inf_q\left\{\lambda\int_{\Theta}d(\theta^*,\theta)\,\mathrm{d}q(\theta)+\frac{D(q\|\pi)}{\beta(n+1)}\right\},
    \end{equation}
where the inf runs over all probability measures $q$ on $\Theta$.
\end{enumerate}
\end{cor}

We end by a short comment on the second statement of Corollary \ref{cor:otbc2}. 
\begin{rem}
\label{rem:futurer}
For any probability measure $q$ on $\Theta$, and since $q\otimes\delta_{\theta^*}$ is the unique coupling between $q$ and $\delta_{\theta^*}$, the upper bound \eqref{eq:otbc2} reads equivalently
\[\esp[\ell(\theta_n,Z)]-\inf_{\theta\in\Theta}\esp[\ell(\theta,Z)]\le\inf_q\left\{\lambda W_1(q,\delta_{\theta^*})+\frac{D(q\|\pi)}{\beta(n+1)}\right\},\]
where $W_1$ denotes the $1$-Wasserstein metric. We believe this form of the upper bound displays an interesting trade-off between the $W_1$ distance to $\delta_{\theta^*}$ and the Kullback-Leibler divergence with respective to the prior $\pi$. Indeed, the upper bound is formally in the same spirit as a proximal gradient step for the functional $D(.\|\pi)$ in the $W_1$ metric. This further suggests to push these investigations in connectiong with the theory of gradient flows in metric spaces and in particular in the space of probability measures \citep{AmbGigSav08}. 
\end{rem}

\subsection{Discussion: On aggregation}
\label{subsec:learn}

 We shortly comment on the use of the previous results in the context of a classical application, the aggregation of predictors. As mentioned earlier, this problem is very well studied in the context of real valued (or euclidean valued) functions. Here, we look at the case of functions taking values in an NPC space.
 
Let $(\mathcal X,d_{\mathcal X})$ be an arbitrary metric space and $(\mathcal Y,d_{\mathcal Y})$ be an NPC space. Let $\{(X_i,Y_i)\}_{i=1}^n$ be i.i.d. random variables with same distribution as (and independent from) a generic pair $(X,Y)$. Consider a set $\Theta$ of functions $\theta:\mathcal X\to\mathcal Y$, such that
\begin{equation}
\label{eq:fsm}
\esp[d^2_{\mathcal Y}(\theta(X),y)]<+\infty,
\end{equation}
for all $y\in\mathcal Y$. Consider the task of building $\theta_n:\mathcal X\to\mathcal Y$ based on the data $\{(X_i,Y_i)\}_{i=1}^n$ such that 
\[\esp[l(Y,\theta_n(X))]-\inf_{\theta\in\Theta}\esp[l(Y,\theta(X))],\]
is as small as possible, for some loss function $l:\mathcal Y^2\to\R$.

First, recall from Example \ref{exm:l2spaces} that, since $(\mathcal Y,d_{\mathcal Y})$ is an NPC space, the set $(\mathcal M,d)$ of all functions $\theta:\mathcal X\to\mathcal Y$ satisfying \eqref{eq:fsm} equipped with metric 
\[d^2(\theta_1,\theta_2)=\esp[d^2_{\mathcal Y}(\theta_1(X),\theta_2(X))],\]
is itself an NPC space. As a result, defining 
\[\ell(\theta,(x,y)):=l(y,\theta(x)),\]
we recover the exact same setting as the one considered in the previous paragraph. 

To apply the previous results, the only detail one may need to check is the fact that convexity assumptions on $l$ translate to convexity properties of $\ell$. This actually holds true given the connection between geodesics in $\mathcal M$ and geodesics in $\mathcal Y$ established in~\citet{Stu03}. Precisely, we have the following lemma.

\begin{lem}
\label{lem:YtoM}
Let $f:\mathcal Y^2\to \R$ be such that, for all $y\in\mathcal Y$, $f(y,.):\mathcal Y\to \R$ is geodesically convex. Then, for all $y\in\mathcal Y$ and $\prob_{X}$-almost all $x\in\mathcal X$, the function $\theta\in\mathcal M\mapsto f(y,\theta(x))$ is geodesically convex.
\end{lem}

\subsection{Discussion: On the estimation of barycenters}
\label{subsec:estimbary}

A second application of interest is the estimation of barycenters. This problem has gained momentum in the past few years and recent contributions are for instance \citet{Sch19,AhiGouPar19,GouParRigStr19} and \citet{Che20}. Despite the surprisingly neat results in these papers, it seems to be still an open question to prove a complexity free rate (no assumption on the covering number of the underlying space) for the estimation of barycenters in NPC spaces with no curvature lower bound (the case of NPC spaces with a curvature lower bound is dealt with in \citealt{GouParRigStr19}). While we do not answer this question here, we provide some comments in this direction.

Note first that the problem of barycenter estimation is closely related to the problem described in paragraph \ref{subsec:otob} with $\Theta=\mathcal M=\mathcal Z$ and $\ell(.,.)=d^2(.,.)$. Indeed, if $\theta^*$ denotes the barycenter of the distribution of $Z$, then Theorem \ref{thm:sturm} implies that, if $\theta_n$ defined as in \eqref{eq:defthetan},  
\begin{align}
\esp[d^2(\theta_n,\theta^*)]&\le \esp[d^2(\theta_n,Z)]-\esp[d^2(\theta^*,Z)]
\nonumber\\
&= \esp[d^2(\theta_n,Z)]-\inf_{\theta\in\Theta}\esp[d^2(\theta,Z)].
\nonumber
\end{align}
In particular, provided $\mathcal M$ has bounded diameter, the second statement of Corollary \ref{cor:otbc2} implies in this case that
\[\esp[d^2(\theta_n,\theta^*)]\le 2\inf_q\left\{\mathrm{diam}(\mathcal M) W_1(q,\delta_{\theta^*})+\mathrm{diam}(\mathcal M)^2\frac{D(q\|\pi)}{n+1}\right\},\]
where the inf runs over all probability measures $q$ on $\Theta$ and where $W_1$ denotes the $1$-Wasserstein metric. 

At this moment, it isn't clear for us if pushing further the analysis of this upper bound could provide a fast rate of order $1/n$ for the EWA based estimator $\theta_n$. This question is left for future research (see Remark \ref{rem:futurer}).

\appendix
\section{Proofs}
\label{sec:proofs}

\subsection{Proof of Theorem \ref{thm:pac}}

For $t\ge 1$, set
\[W_t:=\int_{\Theta}e^{-\beta L_{\theta,t-1}}\,\mathrm{d}\pi(\theta).\]
On the one hand, the Gibbs variational principle implies that
\[\ln W_{T+1}=-\inf_{q}\left\{\beta\int_{\Theta}L_{\theta,T}\,\mathrm{d}q(\theta)+D(q\|\pi)\right\}.\]
On the other hand, since $W_1=1$, we observe that 
\begin{align}
   \ln W_{T+1} &=\sum_{t=1}^T\ln\frac{W_{t+1}}{W_t},
   \nonumber\\
   &=\sum_{t=1}^T\ln\left(\int_{\Theta}e^{-\beta\ell(m_{\theta,t},z_t)}\,\mathrm{d}\pi_t(\theta)\right),
    \nonumber\\
    &\le -\beta \hat L_T,
    \label{eq:weusej}
\end{align}
where the last inequality follows from the concavity of $e^{-\beta\ell(.,z_t)}$, Jensen's inequality and the definition of $\hat m_t$ in \eqref{eq:ewaclassic}. To prove the last statement, it remains to take $q$ to be the Dirac mass $\delta_{\theta^*}$ at any $\theta^*$ minimizing the map $\theta\mapsto L_{\theta,T}$.

 \subsection{Proof of Lemma \ref{lem:lcisec}} 
The proof of the first statement is elementary. It requires to combine simply the fact that the logarithm is increasing and concave. 

\begin{rem}
Before proving the second statement, note that it is elementary in the context of smooth functions on a euclidean space. Indeed, note for example that if $M$ is an interval of the real line, the second derivative of $e^{-\beta f}$ is
\[\beta(\beta (f')^2-f'')e^{-\beta f},\]
which is obviously non-positive whenever $\beta$ satisfies the requirements of the lemma. In the context of geodesic spaces, one cannot refer to usual arguments from differential calculus. The proof we give next is however quite elementary.  
\end{rem}

To prove the second statement, suppose $f:M\to\R$ is geodesically $\alpha$-convex and $L$-Lipchitz for some $\alpha,L>0$. Using the fact that $M$ is complete, it is enough to show that for every $0\le \beta\le ...$ and every geodesic $\gamma:[0,1]\to M$, we have 
\[\frac{1}{2}e^{-\beta f(\gamma(0))}+\frac{1}{2}e^{-\beta f(\gamma(1))}\le e^{-\beta f(\gamma(1/2))},\]
or, equivalently, that 
\begin{equation}
\label{lem:lcisec:e1}
\frac{1}{2}e^{\beta(f(z)-f(x))}+\frac{1}{2}e^{\beta(f(z)-f(y))}\le 1,
\end{equation}
where $x=\gamma(0),y=\gamma(1)$ and $z=\gamma(1/2)$. By geodesic $\alpha$-convexity of $f$, we have that 
\[f(z)\le\frac{1}{2}f(x)+\frac{1}{2}f(y)-\frac{\alpha}{8}d(x,y)^2,\]
so that
\[\beta(f(z)-f(x))\le \frac{\beta}{2}(f(y)-f(x))-\frac{\alpha\beta}{8}d(x,y)^2,\] 
and
\[\beta(f(z)-f(y))\le \frac{\beta}{2}(f(x)-f(y))-\frac{\alpha\beta}{8}d(x,y)^2.\]
For \eqref{lem:lcisec:e1} to be satisfied, it is therefore enough to have 
\begin{equation}
\label{lem:lcisec:e2}
\frac{1}{2}e^{\frac{\beta}{2}(f(y)-f(x))}+\frac{1}{2}e^{\frac{\beta}{2}(f(x)-f(y))}=\cosh(\frac{\beta}{2}(f(x)-f(y)))\le e^{\frac{\alpha\beta}{8}d(x,y)^2}.
\end{equation}
Now, using the fact that $\cosh(u)\le e^{\frac{u^2}{2}}$ for all $u\in \R$, the Lipschitz assumption implies that
\[\cosh(\frac{\beta}{2}(f(x)-f(y)))\le e^{\frac{\beta^2L^2}{8}d(x,y)^2}.\]
As a result, it is enough to have $\beta^2L^2\le \alpha\beta$  i.e. $\beta\le\frac{\alpha}{L^2}$, as required.

\subsection{Proof of Theorem \ref{thm:ewametric1}}
The proof is an immediate adaptation of that of Theorem \ref{thm:pac} using the tools presented in section \ref{sec:geom}. Indeed, letting 
\[W_t:=\int_{\Theta}e^{-\beta L_{\theta,t-1}}\,\mathrm{d}\pi(\theta),\]
for all $t\ge 1$, we can still deduce from the Gibbs variational principle that
\begin{align}
   -\inf_{q}\left\{\beta\int_{\Theta}L_{\theta,T}\,\mathrm{d}q(\theta)+D(q\|\pi)\right\}&=\ln W_{T+1}
   \nonumber\\&=\sum_{t=1}^T\ln\frac{W_{t+1}}{W_t}
   \nonumber\\
   &=\sum_{t=1}^T\ln\left(\int_{\Theta}e^{-\beta\ell(m_{\theta,t},z_t)}\,\mathrm{d}\pi_t(\theta)\right)
      \nonumber\\
   &=\sum_{t=1}^T\ln\left(\int_{\mathcal M}e^{-\beta\ell(x,z_t)}\,\mathrm{d}P_t(x)\right),
    \nonumber
\end{align}
where the last identity follows from the definition of $P_t$ as the pushforward of $\pi_t$ by the map $\theta\in\Theta\mapsto m_{\theta,t}\in\mathcal M$. Finally, using the representation \eqref{eq:ewametric2} of $\hat m_t$ as barycenter of $P_t$, Theorem \ref{thm:jensen} implies that  
\[\sum_{t=1}^T\ln\left(\int_{\mathcal M}e^{-\beta\ell(x,z_t)}\,\mathrm{d}P_t(x)\right)\le -\beta \hat L_T,\]
which completes the proof of the first statement. The last statement follows from the same arguments as in Theorem \ref{thm:pac}.  

\subsection{Proof of Theorem \ref{thm:otbc}}
The proof adapts well known arguments in the context of a euclidean convex set $\mathcal M$. As a matter of fact, the proof only differs from the euclidean setting in that it invokes the less classical Theorem \ref{thm:jensen} to derive the inequality \eqref{eq:otbcmetric}. To ease notation, we identify ($\hat m_1$ with $m_1$ and) $\hat m_i$ with $m_i(Z_1,\dots,Z_{i-1})$ for all $i\ge 2$.

First, since inequality \eqref{eq:boundforotbc} holds for any outcome sequence, we have (almost surely) 
\[\frac{1}{n+1}\sum_{i=1}^{n+1}\ell(\hat m_i,Z_i)-\inf_{\theta\in\Theta}\frac{1}{n+1}\sum_{i=1}^{n+1}\ell(\theta,Z_i)\le \frac{B_{n+1}}{n+1}.\]
Taking the expectation on both sides we deduce that
\begin{align}
    \frac{1}{n+1}\sum_{i=1}^{n+1}\esp[\ell(\hat m_i,Z_i)]&\le\esp[\inf_{\theta\in\Theta}\frac{1}{n+1}\sum_{i=1}^{n+1}\ell(\theta,Z_i)]+\frac{B_{n+1}}{n+1}
    \nonumber\\
    &\le\inf_{\theta\in\Theta}\frac{1}{n+1}\sum_{i=1}^{n+1}\esp[\ell(\theta,Z_i)]+\frac{B_{n+1}}{n+1}
    \nonumber\\
    & =\inf_{\theta\in\Theta}\esp[\ell(\theta,Z)]+\frac{B_{n+1}}{n+1}.
    \label{eq:otbc:1}
\end{align}
Theorem \ref{thm:jensen} implies by definition of $\theta_n$ that, for all $z\in\mathcal Z$, 
\begin{equation}
\label{eq:otbcmetric}
\ell(\theta_n,z)\le \frac{1}{n+1}\sum_{i=1}^{n+1}\ell(\hat m_i,z).
\end{equation}
In particular,
\begin{align}
    \esp[\ell(\theta_n,Z)]&\le \frac{1}{n+1}\sum_{i=1}^{n+1}\esp[\ell(\hat m_i,Z)]
    \nonumber\\
   & = \frac{1}{n+1}\sum_{i=1}^{n+1}\esp[\ell(\hat m_i,Z_i)],
   \label{eq:otbc:2}
\end{align}
where the last identity holds since $\hat m_i$ and $Z_i$ are independent for all $i\ge 1$. The proof then follows from \eqref{eq:otbc:1} and \eqref{eq:otbc:2}.

\subsection{Proof of Theorem \ref{thm:otbc2}}
The proof of the first statement is very similar to the proof of Theorem \ref{thm:otbc} and is avoided for brevity. Formally one should impose an assumption guaranteeing that, for all probability measures $q$ on $\Theta$,
\[\esp[\int_{\Theta}\ell(\theta,Z)\,\mathrm{d}q(\theta)]=\int_{\Theta}\esp[\ell(\theta,Z)]\,\mathrm{d}q(\theta),\]
but we avoid technical discussions on this point since it obviously holds for (say) positive or bounded losses. For the second statement, it is enough to notice that according to the first point, we have for every probability measure $q$ on $\Theta$ that
\begin{align}
    \esp[\ell(\theta_n,Z)]-\esp[\ell(\theta^*,Z)]
    &\le \int_{\Theta}(\esp[\ell(\theta,Z)]-\esp[\ell(\theta^*,Z)])\,\mathrm{d}q(\theta)+\frac{B_{n+1}(q)}{n+1}
    \nonumber\\
     &\le \lambda\int_{\Theta}d(\theta^*,\theta)\,\mathrm{d}q(\theta)+\frac{B_{n+1}(q)}{n+1}.
    \nonumber
\end{align}

\subsection{Proof of Lemma \ref{lem:YtoM}}

The proof follows directly from the following fact. According to~\citet[Proposition 3.10]{Stu03}, a curve $(\Gamma_t)_{t\in[0,1]}$ is a geodesic in $(\mathcal M,d)$ if, and only if, $(\Gamma_t(x))_{t\in[0,1]}$ is a geodesic in $(\mathcal Y,d_{\mathcal Y})$ for $\prob_{X}$-almost all $x\in\mathcal X$.

\bibliography{ewa}
\end{document}